%

\input amstex
\magnification=1200
\loadmsam
\loadmsbm
\loadeufm
\loadeusm
\UseAMSsymbols

\hsize=6.9truein
\hoffset=-0.11truein
\vsize=8.9truein
\voffset=-0.2truein

\def\leftitem#1{\item{\hbox to\parindent{\enspace#1\hfill}}}

\def\boxit#1#2{\hbox{\vrule
	\vtop{%
	\vbox{\hrule\kern#1%
	\hbox{\kern#1#2\kern#1}}%
	\kern#1\hrule}%
	\vrule}}

\def\leaderfill{\leaders\hbox to 1em{\hss.\hss}\hfill}

\parskip=\medskipamount
\document

\input epsf



\centerline{\bf Some Applications of a Multiplicative  
Structure on Simple Loops in Surfaces}

\bigskip

\centerline{Feng Luo}

\bigskip

\centerline{ \it To Joan Birman on her 70th birthday}

\bigskip

\S1. {\bf  Introduction}

Given two transversely intersecting unoriented arcs in an oriented
surface, it is well known that there exists a natural way of resolving
the intersection so that the resolution depends only on the orientation
of the surface and the order of the two arcs (see for instance [BS],
[De], [Wh]). The analogous ``local surgery" operation plays a
prominent role in the recent development of knot invariants (see 
[BL] for instance). 
This local operation  is shown in [Lu] to induce a multiplicative
structure on the space of isotopy classes of curves systems in
surfaces. The goal of this note is to discuss some simple
applications of the multiplicative structure.


The multiplication introduced here is similar to Goldman's Lie bracket
on the free module generated by the conjugacy classes in the 
fundamental group of a surface [Go].  However,  the
multiplication  is  highly non-associative and 
is more geometrically
oriented. We are informed by a referee that in some unpublished work
of Thurston and Penner on train-tracks, they have also discovered
the multiplication.

Suppose $\Sigma = \Sigma_{g,r}$ is a compact oriented surface of genus $g$ with
$r$ boundary components. Following Dehn [De] and Thurston [Th2],
a \it curve system \rm in $\Sigma$ is a proper 1-submanifold so that
each component of the submanifold is homotopically non-trivial
relative to the boundary $\partial \Sigma$. The set $CS(\Sigma)$ of
isotopy classes of curve systems in $\Sigma$ was introduced by Dehn
who called it the \it arithmetic field \rm of the surface. For
instance, if the surface $\Sigma$ is a torus, then the space $CS(\Sigma_{1,0})$ is
naturally isomorphic to the quotient $H_1(\Sigma_{1,0}, \bold Z) -\{0\}/\pm 1$.
Given a lattice $\bold Z^2$, there exists a multiplication on $\bold Z^2  -
\{0\}/
\pm 1$ which is invariant under the sense preserving linear automorphism
group $SL(2, \bold Z)$. Namely, if $\alpha = \pm (x,y)$ and $\alpha'
= \pm (x', y')$, then  $\alpha * \alpha' = \pm (( x,y) + \delta (x', y'))$
where  $\delta$ is the sign of the number $xy' - x'y$ if $xy' - x'y \neq 0$
and $\delta$ is the sign of number $k$ if $(x, y) = k(x', y')$. This
multiplication on $CS(\Sigma_{1,0})$ has a topological interpretation
in terms of the local resolution. In fact, for all oriented surfaces,
there exists a natural multiplication on $CS(\Sigma)$ which is invariant
under the action of the orientation preserving homeomorphisms.

Recall that if $\alpha, \beta \in$ $CS(\Sigma)$, then their \it geometric
intersection number \rm $I(\alpha, \beta)$ is defined to be the
minimal number min$\{ |a \cap b| | a \in \alpha, b \in \beta\}$.
Let $CS_0(\Sigma)$ be the subset of $CS(\Sigma)$ consisting of isotopy classes of
curve systems with no arc components.
The isotopy class of a curve system $c$ is denoted by $[c]$.
We use $\alpha^k$ for $k \in \bold Z_{>0}$ to denote $k$ copies
of $\alpha$.

In [Lu], the following theorem is proved.

{\bf Theorem 1.} \it There exists a  multiplication 
$CS(\Sigma)$$\times$$CS(
\Sigma)$ $\to$ $CS(\Sigma)$ sending
$CS_0(\Sigma)$$\times$ $CS_0(\Sigma)$ to $CS_0(\Sigma)$ and satisfying
the following properties.

(a) It is preserved by the action of the orientation preserving
homeomorphisms.

(b) (Commutative) If $I(\alpha, \beta)$ =0, then $\alpha \beta = \beta \alpha$
and $I(\alpha \beta, \gamma) = I(\alpha, \gamma) + I(\beta, \gamma)$ for
all $\gamma$.
Conversely, if
$\alpha \beta = \beta \alpha$ and $\alpha \in$ $CS_0(\Sigma)$, then $I(\alpha, 
\beta)$=0.

(c) (Associative) If $[c_i] \in$ $CS(\Sigma)$ so that $|c_i \cap c_j| =
I([c_i], [c_j]) $ for
$i,j=1,2,3$,
$i \neq j$, $|c_1 \cap c_2 \cap c_3| =0$,
and there is no contractable region in $\Sigma -(c_1 \cup c_2 \cup c_3)$
bounded either  by three arcs in $c_1$, $c_2$, $c_3$ or four arcs
one in each of $c_1$, $c_2$, $c_3$ and $\partial \Sigma$, 
then $[c_1]([c_2][c_3])= ([c_1][c_2])[c_3]$.

(d) (Cancellation)  If $\alpha \in CS_0(\Sigma)$ and each component of
$\alpha$ intersects $\beta$, then $\alpha(\beta \alpha) = \beta$
and $(\alpha \beta) \alpha = \beta$.
Furthermore, 
$I(\alpha, \alpha \beta)=I
(\alpha, \beta \alpha) = I(\alpha, \beta)$.

(e) For any positive integer $k$, $\alpha ^k \beta^k =(\alpha \beta)^k$.

(f) If $\alpha$ is the isotopy class of a simple closed curve, then the
positive Dehn twist along $\alpha$ sends $\beta$ to $\alpha^k \beta$ where $k =$
 $I(\alpha, \beta)$.

(g) If $\alpha, \beta$ are in $CS_0 (\Sigma)$, then the sum of any
two numbers in $\{I(\alpha, \gamma),  I(\beta, \gamma), I(\alpha \beta, \gamma)
\}$ is at least the third. 
\rm

As one consequence, we prove a convexity result concerning the
intersection number function. It is a combinatorial 
analogous to Kerckhoff's theorem that
the geodesic length function is convex along earthquake path in
the Teichm\"uller space [Ker]. Indeed, it is conceivable that
the multiplication is the extension of Thurston's earthquake
to the measured lamination (see [Bo], and [Pa]). If this holds,
then theorem 2 is also a consequence of Kerckhoff's theorem ([Bo1]).

Given $\alpha$ and $\beta$ in $CS_0(\Sigma)$, 
if $n \in \bold Z_{<0}$, we define
$\alpha^n \beta = \beta \alpha^{-n}$ and $\beta \alpha^n = \alpha^{-n} \beta$.

{\bf Theorem 2} (Convexity). \it  Given $\alpha, \beta \in $ $CS_0(\Sigma)$ and
 $\gamma \in $ $CS(\Sigma)$,
the intersection number function $I(\alpha^n \beta, \gamma)$ is convex
in $n \in \bold Z$. \rm 

{\bf Corollary 3.} \it Suppose $\alpha$ is the isotopy class of a
simple loop in a surface, $\beta \in CS_0(\Sigma)$,
 and $D_{\alpha}$ is the Dehn twist on
$\alpha$. Then the function $I(D^n_{\alpha} (\beta), \gamma)$ is
convex in $n \in \bold Z$. \rm

The paper is organized as follows. In \S2, we  give a short proof
of the main parts of theorem 1. In \S3, we 
prove theorem 2. In \S4, we use the multiplication to
prove several known results.

\it Acknowledgment. \rm I would like to thank F. Bonahon, X.-S. Lin,
Y. Minsky,  C. Series and R. Stong for discussions.  I thank the
referee for his or her nice suggestions.
This work is supported in
part by the NSF.

\S2. {\bf Definitions and Basic Properties of the Multiplication}

Suppose $a$ and $b$ are two unoriented arcs intersecting transversely
at a point 
$p$ in an oriented surface $\Sigma$. Then the \it resolution of 
$a \cup b$ at p from $a$ to $b$ \rm is 
defined as follows. Take any orientation on $a$ and use the
orientation on the surface to determine an orientation
on $b$. Then resolve the intersection according to the orientations.
The resolution is independent of the choice of orientations on $a$.
See  figure 1.

\midspace{0.1cm}
\centerline{\epsfbox{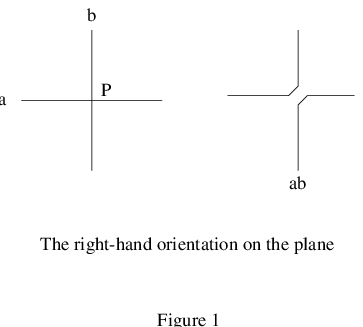}}
\midspace{0.1cm}

Given two curve systems $a, b$ in $\Sigma$ which intersect minimally
within their isotopy classes,  i.e., $I([a], [b]) = |a \cap b|$,
the multiplication $ab$ is defined to be the 1-dimensional submanifold
obtained by resolving all intersection points in $a \cup b$ from 
\it a to b\rm. If $a$ and $b$ are disjoint, then $ab$ is the union
of $a$ with $b$. Our first observation (lemma 4) is that $ab$
is again a curve system. Furthermore, the isotopy class of $ab$
depends only on the isotopy classes of $a$ and $b$.  The multiplication
of two classes $\alpha$ and $\beta$ in $CS(\Sigma)$ is defined to be the
isotopy class of $ab$ where $a \in \alpha$, $b \in \beta$ and
$|a \cap b| = I(\alpha, \beta)$. Note that the positive Dehn
twist $D_{\alpha}$ applied to $\beta$ is the multiplication
$\alpha^k \beta$ where $k = I(\alpha, \beta)$ by the definition. See
figure 2.

\midspace{0.1cm}
\centerline{\epsfbox{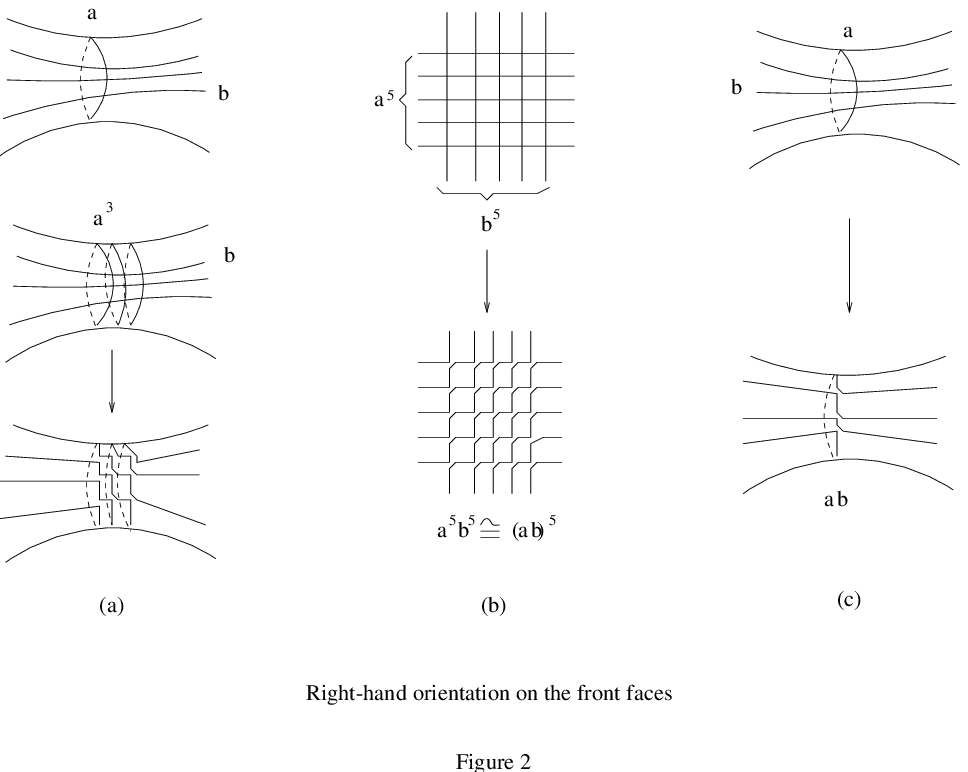}}
\midspace{0.1cm}

{\bf Lemma 4.} \it Suppose $a$ and $b$ are curve systems intersecting
minimally within their isotopy classes.  Then the 1-dimensional
submanifold $ab$ obtained by resolving all intersection points
from $a$ to $b$ is again a curve system. \rm

\it Proof. \rm
We shall prove the lemma in the special case that both $a$ and
$b$ contain no arc components in a closed surface. 
The proof of the general case 
is essentially the same. See [Lu] for details.

Suppose otherwise that the 1-submanifold $ab$ contains a null homotopic
component $c$. We may assume that $c$ is the ``inner-most" component,
i.e., in the interior of the disc  $D$ bounded by $c$, there are no other
components of $ab$.  Let us consider all components of $\Sigma -(a \cup b)$
which are inside $D$, say $A_1,..., A_k$. Each $A_i$ is  an open disc
since $c$ is the inner-most. The boundary of $A_i$ consists of arcs in
$a$ and $b$, and the  \it corners \rm
 of $A_i$ correspond to the intersection points
of $a$ and $b$ in $A_i$. Thus
we may call each $A_i$ a polygon bounded by sides in $a$ and $b$
alternatively. Since $a$ intersects $b$ minimally within their
isotopy classes, each $A_i$ has at least four sides. Now by the definition
of the resolution, the disc $D$ is obtained by resolving corners of
$A_i$'s from $a$ to $b$. Consider the resolutions at the vertices  alone
the boundary of $A_i$. One sees that corners open and
closed alternatively in a  cyclic order on the boundary.
 Form a graph in $D$ by assigning a vertex
in each $A_i$ and joining an edge between two vertices if their corresponding
polygons $A_i$ and $A_j$
have the same corner which is opened by the resolution. Then, on
one hand, the graph is a tree since it is homotopic to the disk.
On the other hand, each vertex of the graph has valency at least two
since the valency of the vertex is half of the number of sides of the
corresponding polygon $A_i$ (by the alternating property).
This contradicts the fact that a tree must have a vertex of valency one.
$\square$

\midspace{0.1cm}
\centerline{\epsfbox{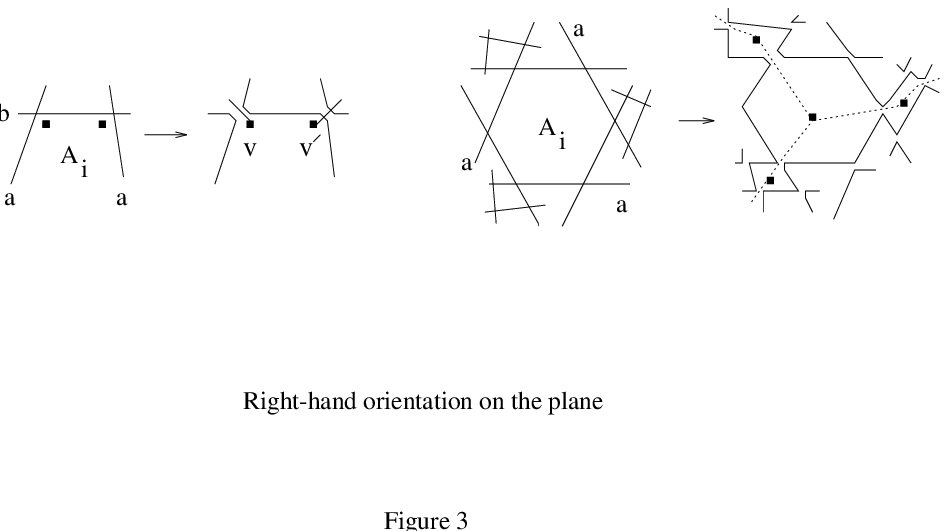}}
\midspace{0.1cm}

Now the  proof of theorem 1 goes as follows. Properties
(a), (e), (f) and the first part of (b) follow from the definition.
Property (g) and the second part of (b) follow from property (d).
Indeed, by the definition of resolution, we have  $I(\alpha  \beta, \gamma)
\leq I(\alpha,  \gamma) + I(\beta, \gamma)$. Now write
$\alpha = \alpha_1 \alpha_2$ with $I(\alpha_1, \alpha_2) = 0$ so
that $I(\alpha_2, \beta) = 0$ and each component of $\alpha_1$
intersects $\beta$. Then by (d), $\beta = (\alpha_1 \beta) \alpha_1$
and $\beta \alpha_2 = \alpha_2 \beta$. Thus, $I(\beta ,\gamma)
=I( (\alpha_1 \beta) \alpha_1, \gamma) \leq I(\alpha_1 \beta, \gamma)
+ I(\alpha_1, \gamma) \leq I(\alpha \beta, \gamma) + I(\alpha, \gamma)$.
To prove part (b),  since $\alpha_1 \alpha_2 \beta = \beta \alpha _1 \alpha_2$, 
we obtain $\alpha_1 \beta = \beta \alpha_1$ by dropping the
component $\alpha_2$. By (d), $\beta = \alpha_1 (\beta \alpha_1)$.
Thus $\beta = \alpha_1^2 \beta$. Now by (d) again 0 = $I(\beta, \beta)
= I(\beta, \alpha^2_1 \beta) = I(\beta, \alpha^2_1) = 2I(\beta, \alpha_1)$.
This shows that $\alpha_1 = \emptyset$ and $I(\alpha, \beta) = I(\beta, \alpha)$.

It remains to prove parts (c) and (d). The key step in proving (c) is
to show that $|c_1 c_2 \cap c_3| = I([c_1 c_2], [c_3])$ and
$|c_1 \cap c_2 c_3| = I([c_1], [c_2 c_3])$ using the non-existence
of triangular and quadrilateral regions in $\Sigma -(\cup_{i=1}^3 c_i)$.
The proof of this is essentially the same as the argument used in
the proof of lemma 4. We refer the reader to [Lu] for details. 
Assuming this, then property (c) follows since both
$([c_1][c_2])[c_3]$ and $[c_1]([c_2][c_3])$ are obtained from 
$c_1 \cup c_2 \cup c_3$ by resolving simultaneously the
intersection points in $ c_i \cap c_j$ from $c_i$ to $c_j$
where $(i,j) = (1,2), (2,3)$ and $(1,3)$. Finally to show property
(d), we note that three classes $\alpha, \beta, \alpha$ satisfy
the non-triangular and non-quadrilateral region  condition in (c). Thus,
by taking $a \in \alpha$, $b \in \beta$ so that $|a \cap b|
= I(\alpha, \beta)$ and $a'$ a parallel copy of $a$, we see that
$I(\alpha, \beta) = I(\alpha, \alpha \beta) = I(\alpha, \beta \alpha)$.
Furthermorem $\alpha (\beta \alpha)$ is represented by
resolving all intersection points in 
$a \cup b \cup a'$
from $a$ to $b$ and $b$ to $a'$. A simple calculation shows
$aba'$ is isotopic to $b$. See figure 4.

\midspace{0.1cm}
\centerline{\epsfbox{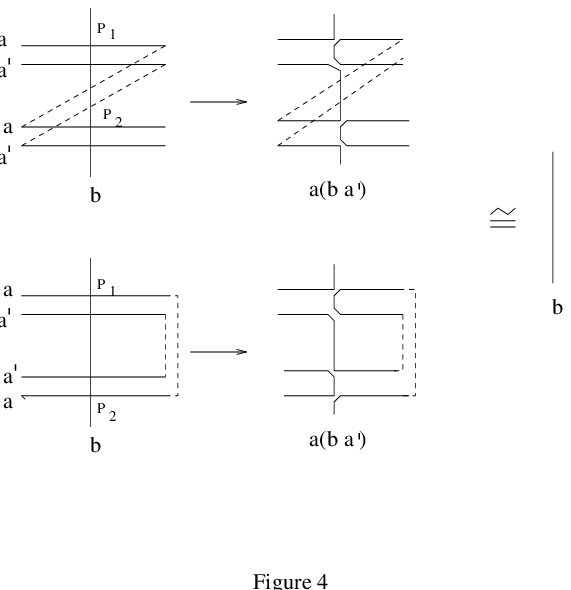}}
\midspace{0.1cm}

\it Remark. \rm The non-triangular region condition has already appeared 
in many papers. The proofs are closely related to the topological
cosine formula in [Bo].

\S3. {\bf Convexity of Intersection Number Functions}

Recall that if $n \in \bold Z_{<0}$, we defined $\alpha^n \beta$ to
be $\beta \alpha^{-n}$. (Note that under this convention, $\alpha^{-2}
\alpha^{5} = \alpha^{7}$.) If each component of $\alpha$ intersects
$\beta$ and $\alpha$ contains no arc component, 
then by theorem  1(d),  we have  $\alpha^n (\alpha^m \beta)
= \alpha^{n+m} \beta$ for all $m,n \in \bold Z$. 

{\bf Theorem 2}.  \it  Given $\alpha, \beta \in $ $CS_0(\Sigma)$ and 
$\gamma \in  $ $CS(\Sigma)$,
the intersection number function $I(\alpha^n \beta, \gamma)$ is convex
in $n \in \bold Z$. \rm

\it Proof. \rm let $f(n) = I(\alpha^n \beta, \gamma)$.
 Since $\alpha^2$, $\beta$ and $\alpha$ satisfy the
non-triangular region condition, by theorem 1(c), we have
$(\alpha^2 \beta) \beta = \alpha ^2 \beta^2 = (\alpha \beta)^2$.
Thus $2 f(1) = 2I(\alpha \beta, \gamma) = I((\alpha \beta)^2, \gamma)$
$ \leq I(\alpha^2 \beta, \gamma) + I(\beta, \gamma)$ $=  f(2) + f(0)$.

Now to prove the convexity of $f$ for general classes, it suffices
to prove it for the case that each component of $\alpha$ intersects
$\beta$.
Indeed, we may write as before $\alpha = \alpha_1 \alpha_2$
with $I(\alpha_1, \alpha_2) =0$ so that $I(\alpha_1, \beta) = 0$ and
each component of $\alpha_2$ intersects $\beta$. Thus $f(n)
= I( \alpha_1 ^{|n|} \alpha_2 ^n \beta, \gamma)$ = $ |n| I(\alpha_1, \gamma)
+ I(\alpha_2^n \beta, \gamma)$ is the sum of two convex functions.

Assume now that each component of $\alpha$ intersects $\beta$. Take
$n_1 < n_2$ in $\bold Z$ so that $n_2 - n_1$ is an even number $2k$.
Let $\alpha' = \alpha^{k}$ and $\beta' = \alpha^{n_1} \beta$. By
the assumption on the components of $\alpha$, we have
$\alpha^{\frac{n_1 + n_2}{2}} \beta = \alpha' \beta'$,
$\alpha^{n_1} \beta = \beta'$, and $\alpha^{n_2} \beta =( \alpha')^2 \beta'$.
Thus $2f(\frac{n_1 + n_2}{2}) \leq f(n_1) + f(n_2)$ follows from 
the first inequality $2f(1) \leq f(0) + f(2)$ for the pair
$(\alpha', \beta')$. $\square$

Corollary 3 follows from theorem 2 and theorem 1(f).

Suppose $a$ is a simple loop in an oriented surface. Then the multiplication
$a^n b$ is the twisting of $b$ $n$ times around $a$.
In particular, given an annulus containing $a$ as a central curve,
if $b_1$ and  $b_2$ are two curve systems in the annulus
so that they have the same end  points, then $b_1$ is isotopic to
$a^n b_2$ by an isotopy which is the identity on the boundary of the
annulus.

This gives a proof of the following.

{\bf Proposition 5.} \it Suppose $\Sigma_{g,n}$ is a compact oriented
surface of genus $g$ with $n$ boundary components $c_1, ..., c_n$. 
Let $\{\alpha_1, ..., \alpha_{3g+n -3}\}$ be a 3-holed sphere decomposition
of the surface. If $\beta_1$ and $\beta_2$ are two curve systems in
$\Sigma_{g,n}$ so that $I(\beta_1, \alpha_i) = I(\beta_2, \alpha_i) >0$ 
and $I(\beta_1, [c_j]) = I(\beta_2, [c_j])$, then there exist integers
$k_1, ..., k_{3g+n-3}$ so that
$$ \beta_1 = \alpha_1 ^{k_1} ... \alpha^{k_{3g+n-3}}_{3g+n-3} \beta_2. \tag *$$

Conversely, if $\beta_1$ and $\beta_2$ are related by the above equation,
then they have the same intersection numbers with $\alpha_i$ and $[c_j]$.
\rm

\it Remarks \rm
1. By theorem 1(c), the above expression $(*)$ is well defined.

2. The result still holds when some $I(\beta_1, \alpha_i) = 0$ by
appropriately interpreting the equation $(*)$.

\it Proof. \rm Take $a_i \in \alpha_i$ and $b_j \in \beta_j$ so that
$a_i \cap a_k = \emptyset$ and $|b_j \cap a_i| = I(\beta_j, \alpha_j)$.
Let $N(a_i)$ be a small regular neighborhood of $a_i$ so that
$N(a_i) \cap N(a_j) = \emptyset$.  Then each component of 
$\Sigma_{g,n} - \cup_{i=1}^{3g+n-3} int(N(a_i))$ is a 3-holed
sphere.  By the classification of the curve systems in the
3-holed sphere, we may assume, after an isotopy of the surface, that
$b_1 = b_2$ in $\Sigma_{g,n} - \cup_{i=1}^{3g+n-3} int(N(a_i))$.
Now apply the above observation on curve systems $b_1 \cap N(a_i)$
and $b_2 \cap N(a_i)$ in the annulus
$N(a_i)$, the result follows. 
$\square$

In terms of the Dehn-Thurston coordinates of $CS(\Sigma_{g,n})$ with
respect to the 3-holed sphere decomposition (see [HP], [FLP]), the
integer $k_i$ in the proposition 
is the difference of the $i$-th twisting coordinate of the
Dehn-Thurston coordinates  of $b_1$ and $b_2$. 

In [LS], the Dehn-Thurston coordinate of the space of measured
laminations is reexamed from the point of view of multiplication.

\S4. {\bf Some Other Applications}

In this section, we give short proofs of three  known
results in surface theory using the multiplications.

{\bf Theorem  6.} \it Let $\alpha$ and $\beta$  be  two isotopy classes
of simple loops in a surface.
\rm

(a) (Ivanov). \it
If the Dehn twists on $\alpha$ and $\beta$ commute,
then $\alpha$ is disjoint from $\beta$, i.e., $I(\alpha, \beta) = 0$. 
\rm

(b) (Thurston). \it If $\alpha$ and $\beta$ fill the 
surface, i.e., $I(\alpha ,\gamma) +
I(\beta, \gamma) > 0$ for all $\gamma \in $ $CS(\Sigma)$, then the
homeomorphism $D_{\alpha}^{-1} D_{\beta}$ has no fixed point
in $CS(\Sigma)$. \rm

Note that  Thurston [Th1] proved a much stronger result that
the homeomorphism is actually pseudo-anosov.

\it Proof. \rm For part (a), let $k =I(\alpha, \beta)$. 
Suppose otherwise that $k >0$.
Consider $D_{\alpha} D_{\beta}(\alpha) = D_{\beta} D_{\alpha}(\alpha)$.
Then we have $\alpha^l (\beta^k \alpha) = \beta^k \alpha$ for some
$l >0$ (indeed, $l = k^2$ by theorem 1(d) that
$I(\alpha, \beta^k \alpha) = I(\alpha, \beta^k) = k^2$). By theorem 1(c),
we obtain, $\alpha^{l-1} \beta^k = \beta^k \alpha$. Left multiply
both sides by $\alpha$ and use the associativity again, we obtain
$\alpha^l \beta^k = \beta^k$. By theorem 1(d),
$0 = I(\beta^k, \beta^k) = I(\alpha^l \beta^k, \beta^k) = I(\alpha^l,
\beta^k) = k^2 l >0$. This is a contradiction. 

For part (b), suppose otherwise that there exists $\gamma \in
$ $CS(\Sigma)$ so that $D_{\alpha}^{-1} D_{\beta}(\gamma) = \gamma$. Then
$\alpha^n \gamma = \beta^m \gamma$ for some $m,n \in \bold Z_{\geq 0}$.
Multiply from the left by $\gamma$. By theorem 1(d), we obtain
$\alpha^n \gamma' = \beta^m \gamma''$ where $I(\gamma', \alpha) =0$
and $I(\beta, \gamma'') =0$.  One of the classes $\gamma'$ and
$\gamma''$ is non-empty since $\alpha$, $\beta$  are surface filling.
Say $\gamma' \neq 0$. Then the equation $\alpha^n \gamma' = \beta^m \gamma''$
shows that $I(\gamma', \beta)
=0$ which contradicts the  surface filling assumption.
$\square$

The following result was obtained in [FLP] (proposition 1  in the
appendix of expos\'e 4) and [MS] (proposition III3.4).

{\bf Proposition 7} ([FLP], [MS]). \it Suppose  $\alpha_1, ..., \alpha_m$
are isotopy classes of  simple loops so that $I(\alpha_i, \alpha_j) =0$
for all $i,j$. Then  for all $\beta \in CS_0(\Sigma)$ and $\gamma \in
CS(\Sigma)$, we have
$$ \sum_{i=1}^m I(\alpha_i, \beta) I(\alpha_i, \gamma) - I(\beta, \gamma)
\leq I(D_{\alpha_1} ... D_{\alpha_m} \beta, \gamma)
\leq  \sum_{i=1}^m I(\alpha_i, \beta) I(\alpha_i, \gamma) + I(\beta, \gamma).
$$
\rm

\it Proof. \rm Let $\alpha = \alpha_1^{k_1}....\alpha_m^{k_m}$ where
$k_i = I(\alpha_i, \beta)$. Then $I(\alpha, \gamma) 
= \sum_{i=1}^m I(\alpha_i, \beta) I(\alpha_i, \gamma)$. The result is
a consequence of theorem 1(g).
$\square$

\bigskip

\centerline{\bf References}

\bigskip

[Bo] Bonahon, F.: Earthquakes on Riemann surfaces and on 
measured  geodesic laminations. Trans. Amer. Math. Soc. 330 (1992), no. 1, 69--95. 

[Bo1] Bonahon, F.: Private communications.

[BL] Birman, J.,  Lin X.-S.:
Knot  polynomials and Vassiliev's invariants. Invent.
Math. 111 (1993), no. 2, 225--270. 

[BS] Birman, J., C. Series:
 An algorithm for simple curves on surfaces. J. London
Math. Soc. (2) 29 (1984), no. 2, 331--342. 

[De] Dehn, M.: Papers on group theory and topology. J. Stillwell (eds.).
 Springer-Verlag, Berlin-New York, 1987.

[FLP] Fathi, A., Laudenbach, F., Poenaru, V.: Travaux de Thurston sur les
surfaces. Ast\'erisque  66-67, Soci\'et\'e Math\'ematique de France, 1979.

[Go] Goldman, W.: 
Invariant functions on Lie groups and Hamiltonian flows of surface
group representations. Invent. Math. 85 (1986), no. 2, 263--302.

[Iv] Ivanov, N.:
Automorphisms of Teichm\"uller modular groups. Topology and
geometry---Rohlin Seminar, 199--270, Lecture Notes in Math., 1346, Springer, Berlin-New York, 1988. 

[Ker] Kerckhoff, S.: The Nielsen realization problem. Ann. of
Math.  117 (1983), no. 2, 235--265.

[Lu] Luo, F.: Simple loops and their intersections on surfaces, preprint, 1997.

[LS] Luo, F., Stong R.: Dehn-Thurston coordinate for curve systems in surfaces,
in preparation.

[MS] Morgan, J., Shalen, P.: Valuations, trees, and degenerations of
hyperbolic structures, I. Ann. of Math., 120 (1984), 401-476.

[Pa] Papadopoulos, A.: On Thurston's boundary of Teichm\"uller
space and the extension of earthquakes. Topology Appl. 41 (1991), no. 3,
147--177.

[PH] Penner, R.,  Harer, J.: Combinatorics of train tracks.
Annals of Mathematics Studies, 125. Princeton University Press, Princeton,
NJ, 1992.

[Th1] Thurston, W.: On the geometry and dynamics of diffeomorphisms of
surfaces. Bul. Amer. Math. Soc. 19 (1988) no 2, 417-438.       

[Wh] Whitney, H.:
Collected papers.  Edited and with a preface by James Eells and
Domingo Toledo. Contemporary Mathematicians. Birkhauser Boston, Inc., Boston, MA, 1992. 

\bigskip

Department of Mathematics

Rutgers University

Piscataway, NJ 08854

email: fluo\@math.rutgers.edu

\end